\input amstex
\documentstyle{amsppt}
\document
\topmatter
\title
Compact holomorphically pseudosymmetric K\"ahler manifolds.
\endtitle
\author
W{\l}odzimierz Jelonek
\endauthor

\abstract{The aim of this paper is to present the first  examples
of compact, simply connected  holomorphically pseudosymmetric
K\"ahler manifolds. }
\thanks{MS Classification: 53C55,53C25. Key words and phrases:
K\"ahler manifold, holomorphically pseudosymmetric curvature
tensor}\endthanks
 \endabstract
\endtopmatter

\define\DE{\Cal D^{\perp}}

\define\n{\nabla}
\define\om{\omega}

\define\k{\diamondsuit}
\define\th{\theta}

\define\g{\gamma}
\define\lb{\lambda}

\define\1{D_{\lb}}
\define\2{D_{\mu}}
\define\0{\Omega}

\define\De{\Cal D}

\define\m{(M,g,J)}
\define \E{\Cal E}

{\bf 1. Introduction.}  The holomorphically pseudosymmetric
K\"ahler manifolds  were defined by Z. Olszak in [O-1] in 1989 and
studied in [D], [H], [Y]. A K\"ahler manifold $\m$ is called
holomorphically pseudosymmetric if its curvature tensor $R$
satisfies the condition
$$R.R=f\Pi.R,$$
where $R,\Pi$ act as  derivations of the tensor algebra,
$\Pi(U,V)X=\frac14(g(V,X)U-g(U,X)V+g(JV,X)JU-g(JU,X)JV-2g(JU,V)JX)$
is the K\"ahler type curvature tensor of constant holomorphic
sectional curvature and $f\in C^{\infty}(M)$ is a smooth function.
A Riemannian manifold $(M,g)$ is called semisymmetric if $R.R=0$.
Until now  examples of compact, holomorphically pseudosymmetric
and not semisymmetric, K\"ahler manifolds were not known. Z.
Olszak in [O-2] proved that if $\m$ is compact, holomorphically
pseudosymmetric, has constant scalar curvature $\tau$ and $f\ge 0$
then $M$ is locally symmetric - $\n R=0$, where $\n$ is the
Levi-Civita connection of $\m$, and thus semisymmetric. He also
constructed many examples of non-compact holomorphically
pseudosymmetric K\"ahler manifolds, among others holomorphically
pseudosymmetric, not semisymmetric,  K\"ahler-Einstein manifolds
with $f>0$. The QCH K\"ahler manifolds are the K\"ahler manifolds
admitting a smooth, two-dimensional, $J$-invariant distribution
$\De$ whose holomorphic curvature $K(\pi)=R(X,JX,JX,X)$ of any
$J$-invariant $2$-plane $\pi\subset T_xM$, where $X\in \pi$ and
$g(X,X)=1$, depends only on the point $x$ and the number
$|X_{\De}|=\sqrt{g(X_{\De},X_{\De})}$, where $X_{\De}$ is the
orthogonal projection of $X$ on $\De$. In this case  we have
$$R(X,JX,JX,X)=\phi(x,|X_{\De}|)$$ where $\phi(x,t)=a(x)+b(x)t^2+c(x)t^4$ and
 $a,b,c$ are smooth functions on $M$. Also $R=a\Pi+b\Phi+c\Psi$
 for certain curvature tensors $\Pi,\Phi,\Psi\in \bigotimes^4\frak X^*(M)$
  of K\"ahler type. The investigation of  QCH K\"ahler manifolds,  was
started by G. Ganchev and V. Mihova in [G-M-1],[G-M-2]. Compact,
simply connected QCH K\"ahler manifolds were partially classified
by the author  in [J]. We shall show in the present paper that any
QCH K\"ahler manifold is holomorphically pseudosymmetric. In that
way using the classification result from [J] we obtain many
compact, simply connected examples of holomorphically
pseudosymmetric K\"ahler manifolds which are not semisymmetric.

{\bf 2. The curvature tensor of a QCH K\"ahler manifold.} Let us
denote by $\E$ the distribution $\DE$, which is a
$2(n-1)$-dimensional, $J$-invariant distribution. By $h$ we shall
denote the tensor $h=g\circ (p_{\De}\times p_{\De})$, where
$p_{\De}$ are the orthogonal projections on $\De$.  By
$\0=g(J\cdot,\cdot)$ we shall denote the K\"ahler form of $\m$, by
$\om$  the K\"ahler form of $\De$ i.e. $\om(X,Y)=h(JX,Y)$. We
shall recall some results proved by Ganchev and Mihova in [G-M-1].
Let $R(X,Y)Z=([\n_X,\n_Y]-\n_{[X,Y]})Z$ and let us write
$$R(X,Y,Z,W)=g(R(X,Y)Z,W).$$ We shall identify $(1,3)$ tensors
with $(0,4)$ tensors in this way. If  $R$ is the curvature tensor
of a QCH K\"ahler manifold $\m$, then
$$R=a\Pi+b\Phi+c\Psi,\tag 2.1$$
where  $a,b,c\in C^{\infty}(M)$ and    $\Pi$ is the standard
K\"ahler tensor of constant holomorphic curvature
$$\gather \Pi(X,Y,Z,U)=\frac14(g(Y,Z)g(X,U)-g(X,Z)g(Y,U)\tag 2.2\\+g(JY,Z)g(JX,U)-g(JX,Z)g(JY,U)-2g(JX,Y)g(JZ,U)),\endgather $$
the tensor $\Phi$ is as follows
$$\gather \Phi(X,Y,Z,U)=\frac18(g(Y,Z)h(X,U)-g(X,Z)h(Y,U)\tag 2.3\\+g(X,U)h(Y,Z)-g(Y,U)h(X,Z)
+g(JY,Z)\om(X,U)\\-g(JX,Z)\om(Y,U)+g(JX,U)\om(Y,Z)-g(JY,U)\om(X,Z)\\
-2g(JX,Y)\om(Z,U)-2g(JZ,U)\om(X,Y)),\endgather$$ and
$$\Psi(X,Y,Z,U)=-\om(X,Y)\om(Z,U)=-(\om\otimes\om)(X,Y,Z,U).\tag 2.4$$

 Let $V=
 (V,g,J)$ be a  real $n$ dimensional, $n=2k$,  vector space with a
 scalar product $g$ and a
complex structure $J$ such that $g(J\cdot,J\cdot)=g(\cdot,\cdot)$.
Let  $D$ be a 2-dimensional, $J$-invariant subspace  of $V$. By
 $E$ we denote its orthogonal complement in $V$. Then $V=D\oplus
 E$.
 The tensors $\Pi,\Phi,\Psi$ given above are of K\"ahler
type.  We have the following relations for a unit vector $X\in V$:
$\Pi(X,JX,JX,X)=1,\Phi(X,JX,JX,X)=|X_{D}|^2,\Psi(X,JX,JX,X)=|X_{D}|^4$,
where $X_D$ means the orthogonal projection of a vector $X$ on the
subspace $D$ and $|X|=\sqrt{g(X,X)}$. Consequently for a tensor
$(2.1)$ defined on $V$ we have
$$R(X,JX,JX,X)=\phi(|X_D|)$$ where $\phi(t)=a+bt^2+ct^4$.

\medskip
{\bf 3. QCH K\"ahler manifolds are holomorphically
pseudosymmetric.} Let us recall that Z. I. Szab\'o [Sz-1], [Sz-2]
classified semisymmetric Riemannian manifolds, i.e.,  manifolds
whose curvature tensor $R$ satisfies the condition $R.R=0$.  In
particular it turns out that the products
$(\Bbb{CP}^n,g_{can})\times (\Sigma,g)$ where
$(\Bbb{CP}^n,g_{can})$ is the complex projective space with the
standard Fubini-Study metric and $(\Sigma,g)$ is any Riemannian
surface are semisymmetric K\"ahler manifolds. It is a trivial fact
that the product $(M,g)= (M(k),g_k)\times (\Sigma(l),g_l)$, where
$(M(k),g_k)$ is a $(n-2)$ space, $n\ge 4$ of constant holomorphic
sectional curvature $k$ and $(\Sigma(l),g_l)$ is a 2-dimensional
Riemannian surface of constant sectional curvature $l$, is a
locally symmetric (and, in particular semisymmetric) K\"ahler
manifold. Using this fact we shall prove
\medskip
{\bf Theorem 1.  } {\it Any QCH K\"ahler manifold is
holomorphically pseudosymmetric, more precisely if
$R=a\Pi+b\Phi+c\Psi$ is a curvature tensor of a QCH K\"ahler
manifold $\m$ then}

$$ R.R=(a+\frac b2)\Pi.R.\tag 3.1$$

\medskip
{\it Proof.}  We  shall first  prove that the tensors
$\Pi,\Phi,\Psi$ satisfy the following relations:
$$ \gather  2\Phi.\Phi=\Phi.\Pi+\Pi.\Phi,  \Psi.\Psi=0\tag 3.2\\
\Psi.\Pi+\Pi.\Psi=2(\Phi.\Psi+\Psi.\Phi).\endgather$$

Let $V=T_{x_0,y_0}M$ for $(M,g)= (M(k),g_k)\times (\Sigma(l),g_l)$
where $(M(k),g_k)$ is a $(n-2)$ space, $n\ge 4$, of constant
holomorphic sectional curvature $k$ and $(\Sigma(l),g_l)$ is a
2-dimensional Riemannian surface of constant sectional curvature
$l$, $x_0\in M(k),y_0\in \Sigma(l)$. If
$E=T_{x_0}M,D=T_{y_0}\Sigma$ then $V=E\oplus D$ is an orthogonal
sum. Let $R,R^k,R^l$ be the curvature tensors of $(M,g),
(M(k),g_k),(\Sigma(l),g_l)$ respectively. If $X=Y+Z$, where $X\in
V,Y\in E, Z\in D$, is a unit vector,  then
$$R(X,JX,JX,X)=R^k(Y,JY,JY,Y)+R^l(Z,JZ,JZ,Z)=(1-|X_{D}|)^2k+|X_{D}|^4  l,     $$
where $X_{D}=Z$ is the orthogonal  projection of $X$ on $D$. Thus

$$R(X,JX,JX,X)=k-2k|X_{D}|^2+(l+k)|X_{D}|^4.$$
Consequently
$$R=k\Pi-2k\Phi+(l+k)\Psi.\tag 3.3$$

Let us take $l=-k$.  Then   $R=k(\Pi-2\Phi)$. Since $R.R=0$ we
obtain
$$ (\Pi-2\Phi).(\Pi-2\Phi)=0$$ and consequently, since $\Pi.\Pi=0$
$$2\Phi.\Phi=\Phi.\Pi+\Pi.\Phi.\tag 3.4$$

Now let $k=1,d=l+k\ne 0$.  Then $R=\Pi-2\Phi+d\Psi$.  Since
$R.R=0$ we get $(\Pi-2\Phi+d\Psi).(\Pi-2\Phi+d\Psi)=0$ and
$$-2\Pi.\Phi+4\Phi.\Phi+d\Pi.\Psi-2\Phi.\Pi-2d\Phi.\Psi+d\Psi.\Pi-2d\Psi.\Phi+d^2\Psi.\Psi=0$$

One can easily check that  $\Psi.\Psi=0$.  Hence
$\Psi.\Pi+\Pi.\Psi=2(\Phi.\Psi+\Psi.\Phi)$.  Note that these
formulas obtained for particular K\"ahler manifolds are purely
algebraic in $g,h,\om,\0$ and thus remain true for general tensors
$\Pi,\Phi,\Psi$ on a complex vector space $(V,J)$, such that $V
=E\oplus D, JE=E,JD=D, \dim D=2$  and a $J$-invariant  metric $g$
on $V$ such that $g(E,D)=0$ and with $h=g_{|D},\om=h(J.,.)$.
\medskip
Now let $\m$ be a QCH K\"ahler manifold with curvature tensor
$R=a\Pi+b\Phi+c\Psi$.  Then
$$\gather R.R=(a\Pi+b\Phi+c\Psi).(a\Pi+b\Phi+c\Psi)=a\Pi.R+
ab\Phi.\Pi+b^2\Phi.\Phi\tag 3.5\\+bc\Phi.\Psi
+ac\Psi.\Pi+bc\Psi.\Phi
=a\Pi.R+\frac{b^2}2(\Phi.\Pi+\Pi.\Phi)+bc(\Phi.\Psi+\Psi.\Phi)\\+ab\Phi.\Pi+ac\Psi.\Pi=
a\Pi.R+\frac{b^2}2(\Phi.\Pi+\Pi.\Phi)+\frac{bc}2(\Psi.\Pi+\Pi.\Psi)+a((b\Phi+c\Psi).\Pi)\\
=a\Pi.R+\frac b2\Pi.(b\Phi+c\Psi)+\frac b2 (b\Phi+c\Psi).\Pi+a
R.\Pi\\=(a+\frac b2)\Pi.R+(a+\frac b2)R.\Pi=(a+\frac
b2)\Pi.R,\endgather$$ since $R.\Pi=0$ in view of the equality $\n
\Pi=0$.

We also give another, algebraic proof of relations (3.2). Note
that $\Psi(X,Y)Z=-\om(X,Y)J\circ\pi Z$ where $\pi:V\rightarrow D$
is an orthogonal projection and $J$ is a complex structure. Note
that $\pi\circ J=J\circ \pi$.  It is easy to see that
$\Psi.\Pi=\Psi.\Phi=\Psi.\Psi=0$. We shall show the formula
$$\Pi.\Psi=2\Phi.\Psi.$$

We have $\Psi=-\om\otimes\om$ and consequently
$$\Pi.\Psi=-\Pi.\om\otimes\om-\om\otimes\Pi.\om.$$

Thus it is enough to show that $\Pi.\om=2\Phi.\om$. Note that
$$\gather\Pi(U,V)W=\frac14(g(V,W)U-g(U,W)V+g(JV,W)JU-g(JU,W)JV-\\2g(JU,V)JW)\endgather$$
and consequently
$$\gather-4\Pi(U,V).\om(X,Y)=\om(U,Y)g(V,X)+g(V,Y)\om(X,U)\tag 3.6\\-g(U,X)\om(V,Y)-g(U,Y)\om(X,V)
+g(JV,X)\om(JU,Y)+
g(JV,Y)\om(X,JU)\\-g(JU,X)\om(JV,Y)-g(JU,Y)\om(X,JV).\endgather$$

Note that $$\gather 8\Phi(U,V)W=g(V,W)\pi U-g(U,W)\pi V+h(V,W)
U-h(U,W) V\\+g(JV,W)\pi\circ J U-g(JU,W)\pi\circ J
V+\om(V,W)JU\\-\om(U,W)JV-2g(JU,V)\pi\circ J
W-2\om(U,V)JW.\endgather$$

thus
$$\gather
-8\Phi(U,V).\om(X,Y)=g(V,X)\om(U,Y)+g(V,Y)\om(X,U)\tag 3.7
\\-g(U,X)\om(V,Y)-g(U,Y)\om(X,V)
+h(V,X)\om(U,Y)+h(V,Y)\om(X,U)\\-h(U,X)\om(V,Y)-h(U,Y)\om(X,V)
+g(JV,X)\om(JU,Y)+\\g(JV,Y)\om(X,JU)-g(JU,X)\om(JV,Y)-g(JU,Y)\om(X,JV)
\\+\om(V,X)\om(JU,Y)+\om(V,Y)\om(X,JU)-\om(U,X)\om(JV,Y)\\-\om(U,Y)\om(X,JV)
=g(V,X)\om(U,Y)+g(V,Y)\om(X,U)\\-g(U,X)\om(V,Y)-g(U,Y)\om(X,V)
+g(JV,X)\om(JU,Y)\\+g(JV,Y)\om(X,JU)-g(JU,X)\om(JV,Y)-g(JU,Y)\om(X,JV).\endgather$$
Hence  $$\Pi(U,V).\om(X,Y)=2\Phi(U,V).\om(X,Y).$$ It is easy to
see that $\Pi(U,V).h(X,Y)=2\Phi(U,V).h(X,Y)$ and by easy
calculations we get
$$\Phi.g=0,\Phi.J=0,\Phi.\0=0,$$ and consequently $\Phi.\Pi=0$. Also by
direct calculations one can see that $\Pi.g=0,\Pi.J=0,\Pi.\0=0$.
Since $(\Pi-2\Phi).\om=0$ and $(\Pi-2\Phi).h=0$ it follows that
$$(\Pi-2\Phi).\Phi=0.$$ Thus $\Pi.\Phi=2\Phi.\Phi$. Summarizing we
have:
$$\gather
\Pi.\Pi=\Phi.\Pi=\Psi.\Pi=\Psi.\Phi=\Psi.\Psi=0,\tag 3.8\\
\Pi.\Phi=2\Phi.\Phi,\Pi.\Psi=2\Phi.\Psi\endgather$$$\k$

Now it is easy to prove
\medskip
{\bf Theorem 2. }  {\it There exist infinitely many mutually
nonhomothetic holomorphically pseudosymmetric compact, simply
connected K\"ahler manifolds which are not semisymmetric. The
class of compact simply connected QCH K\"ahler manifolds is
included in the class of essentially holomorphically
pseudosymmetric K\"ahler manifolds.}
\medskip
{\it Proof.}  Note that in the case of a QCH K\"ahler manifold
$\m$ constructed in [J] we have
$M_k=P_k\times_{S^1}\Bbb{CP}^1,k\in \Bbb{N}$, where
$p:P_k\rightarrow \Bbb{CP}^{n-1}$ is the $S^1$-bundle over
$\Bbb{CP}^{n-1}$ with a connection form $\th_k$ and $n\ge 2$.
Hence $M$ is a holomorphic $\Bbb{CP}^1$-bundle over
$\Bbb{CP}^{n-1}$. An open and dense subset of $M_k$ is isometric
to the manifold the manifold $(0,L)\times P_k$ with the metric
$$g_k=dt^2+f(t)^2\th_k^2+r(t)^2p^*h,\tag 3.9$$ where $f=\frac{2rr'}s$,
$h$ is the metric of constant holomorphic sectional curvature on
$\Bbb{CP}^{n-1}, s=\frac {2k}n, k\in \Bbb{N}$. The metric (3.9)
extends on the whole of $M_k$ if the function $r$ is positive and
smooth on $(0,L)$, even at the points $0,L$, with $r'>0$ on
$(0,L)$, $r'(0)=r'(L)=0$ and
$$2r(0)r''(0)=s,2r(L)r''(L)=-s.\tag 3.10$$
There exist uncountably many functions $r$ satisfying these
conditions. The curvature tensor of $(M_k,g_k,J)$ is of the form
$R=a\Pi+b\Phi+c\Psi$ and

$$a+\frac b2=4(\frac{(r')^2}{r^2}-\frac{f'r'}{rf})=-4\frac{r''}r.
\tag 3.11$$ Since $r$ satisfies (3.10) the function $r''$ does not
have constant sign and it is clear that the function $a+\frac b2$
changes sign on $M$. From the results of Szab\'o ([Sz-2]) it is
clear that $(M_k,g_k,J)$ is not semisymmetric. One can also
directly verify that $\Pi.R\ne 0$ for $(M_k,g_k,J)$ and hence
$R.R\ne 0$. $\k$

\bigskip

{\it Acknowledgments.} The author would like to thank Professor
Andrzej Derdzi\'nski for his valuable remarks which improved the
paper.

\bigskip
\centerline{\bf References.}

\medskip
[D] R. Deszcz ,  {\it On pseudosymmetric spaces},  Bull. Soc.
Math. Belg., S\'er.A, (1992),  44, 1-34.
\medskip
[H] Hotlo\'s M.,  {\it On holomorphically pseudosymmetric
K\"ahlerian manifolds}, In: Geometry and topology of submanifolds,
VII (Leuven 1994,Brussels 1994), 139-142, World Sci. Publ. River
Edge, NJ, (1995)
\par
\medskip
\cite{G-M-1} G.Ganchev, V. Mihova {\it K\"ahler manifolds of
quasi-constant holomorphic sectional curvatures}, Cent. Eur. J.
Math. 6(1),(2008), 43-75.
\par
\medskip
\cite{G-M-2} G.Ganchev, V. Mihova {\it Warped product K\"ahler
manifolds and Bochner-K\"ahler metrics}, J. Geom. Phys. 58(2008),
803-824.
\par
\medskip
 [J]  W.  Jelonek , {\it K\"ahler manifolds with
quasi-constant holomorphic curvature}, Ann. Global Analysis and
Geom., (2009), DOI : 10.1007/s10455-009-9154-z.

\medskip
[O-1] Z.  Olszak ,  {\it Bochner flat K\"ahlerian manifolds with a
certain condition on the Ricci tensor}, Simon Stevin, (1989),  63,
295-303.

\medskip
[O-2] Z.  Olszak ,  {\it On compact holomorphically
pseudosymmetric K\"ahlerian  manifolds}, preprint (2009).
\medskip
[Sz-1] Z.I. Szab\'o,  {\it Structure theorems on Riemannian spaces
satisfying $R(X,Y).R$ $=0$. I. The local version}, J. Diff. Geom.,
(1982), 17, 531-582.

\medskip
[Sz-2] Z.I. Szab\'o,  {\it Structure theorems on Riemannian spaces
satisfying $R(X,Y).R$ $=0$. II. Global versions},  Geom. Dedicata,
(1985), 19, 65-108.
\medskip
[Y]  Yaprak S., {\it Pseudosymmetry type curvature conditions on
K\"ahler hypersurfaces}, Math. J.  Toyama Univ., (1995), 18,
107-136.

\medskip Institute of Mathematics

Technical University of Cracow

 Warszawska 24

31-155 Krak\'ow,POLAND.

E-mail address: wjelon\@pk.edu.pl

\end